\documentclass[12pt]{amsart}
\usepackage{amscd,amsmath,amssymb,amsfonts}

\theoremstyle{plain}
\newtheorem{thm}{Theorem}
\newtheorem{lem}[thm]{Lemma}

\newtheorem{prop}[thm]{Proposition}

\newtheorem{defn}[thm]{Definition}

\numberwithin{thm}{section}
\numberwithin{equation}{section}

\newcommand{\sA}{{\mathcal A}}
\newcommand{\sB}{{\mathcal B}}

\newcommand{\sD}{{\mathcal D}}
\newcommand{\sE}{{\mathcal E}}

\newcommand{\sH}{{\mathcal H}}

\newcommand{\sK}{{\mathcal K}}
\newcommand{\sL}{{\mathcal L}}

\newcommand{\sO}{{\mathcal O}}
\newcommand{\sP}{{\mathcal P}}
\newcommand{\sQ}{{\mathcal Q}}

\newcommand{\Z}{{\mathbb Z}}

\begin{document}

\title[Determinant bundle]{Determinant bundle in a family of
curves, after A. Beilinson and V. Schechtman} 
\author{H\'el\`ene Esnault}
\address{Mathematik,
Universit\"at Essen, FB6, Mathematik, 45117 Essen, Germany}
\email{esnault@uni-essen.de}
\author{I-Hsun Tsai}
\address{Department of Mathematics, National Taiwan University, Taipei, Taiwan}
\email{ihtsai@math.ntu.edu.tw}  
\date{October 7, 1999}
\begin{abstract}
Let $\pi: X \to S$ be a smooth projective family of curves over a smooth
base $S$ over a field of characteristic 0, together with a
bundle $E$ on X. Then A. Beilinson and
V. Schechtman define in \cite{BS} a beautiful ``trace complex ''
$^{tr}\sA^{\bullet}_E$ on $X$, the 0-th relative cohomology of
which describes the Atiyah algebra of the determinant bundle of $E$ on
S. Their proof reduces the general case to the acyclic one. In
particular, one needs a comparison of
$R\pi_*(^{tr}\sA^{\bullet}_F)$ for $F= E$ and $F=E(D)$ where $D$
is \'etale over $S$ (see Theorem 2.3.1, reduction ii) in
\cite{BS}). In this note, we analyze this reduction in more details
and correct a point.
\end{abstract}
\maketitle

\section{Introduction}
Let $\pi:X\to S$ be   
a smooth projective morphism of relative dimension 1 
over a smooth base $S$ over a field $k$ of characteristic 0.   
One denotes by $T_X$ and the tangent sheaf over $k$, by $T_{X/S}$  
the relative tangent sheaf, and by $\omega_{X/S}$ 
the relative dualizing sheaf.     
For  an algebraic vector bundle $E$ on $X$,
one writes $E^{\circ}=E^*\otimes_{\sO_X}\omega_{X/S}$.
Let ${\sD}iff(E,E)$ (resp. ${\sD}iff({E/S,E/S}) \subset {\sD}iff(E,E)$) 
be the sheaf of first order
(resp. relative) differential  operators on $E$ and
$\epsilon:{\sD}iff(E,E)\to {\sE}nd(E)\otimes_{{\sO}_X}T_X$
be the symbol map.  The Atiyah algebra  
$\sA_E:=\{a\in {\sD}iff(E,E)|\,
\epsilon(a)\in {\rm id}_E\otimes_{\sO_X}T_X\}$ of $E$ is the subalgebra
of ${\sD}iff(E,E)$ consisting of the differential operators 
for which the symbolic
part is a homothety. Similarly  
the relative Atiyah algebra 
$\sA_{E/S}\subset \sA_E$  of $E$ consists of those differential
operators with symbol
in ${\rm id}_E\otimes_{\sO_X} T_{X/S}$ and
$\sA_{E,\pi}\subset \sA_E$ with symbols in 
$T_{\pi}=d\pi^{-1}(\pi^{-1}T_S)\subset T_X$.     
Let $\Delta \subset X\times_S X$
be the diagonal.   
Then there is a canonical
sheaf isomorphism ${\sD}iff(E/S,E/S)\cong \frac{E\boxtimes E^{\circ}(2\Delta)}
{E\boxtimes E^{\circ}}$ (see \cite{BS}, section 2)
which is locally written as follows. 
Let $x$ be a local coordinate of $X$ at a point $p$, and $(x,y)$ be the induced
local coordinates on $X \times _S X$ at $(p,p)$, such that
the equation of $\Delta$ becomes $x=y$. Let $e_i$ be a local
basis of $E$, $e_j^*$ be its local dual basis. Then the action
of $$P= \sum_{i,j} e_i \otimes e_j^* \frac{P_{ij}(x,y)}{(x-y)^2} dy$$
on $s=\sum_\ell s_\ell(y) e_\ell$ is 
\begin{gather}
P(s)= \sum_i e_i \sum_j (P_{ij}^{(1)}(x,0)s_j(x) + 
P_{ij}(x,x)s_j ^{(1)}(x, 0)),
\end{gather}
where 
\begin{gather} 
P_{ij}(x,y)= P_{ij}(x,x) + (y-x)P_{ij}^{(1)}(x, y-x) \notag \\
 s_j(y)=s_j(x) + (y-x)s_j ^{(1)}(x, y-x). \notag
\end{gather}

Beginning with 
\begin{gather}
0\to \frac{E\boxtimes E^{\circ}}{E\boxtimes E^{\circ}(-\Delta)}
\to
\frac{E\boxtimes E^{\circ}(2\Delta)}{E\boxtimes E^{\circ}(-\Delta)}
\to
{\sD}iff(E/S,E/S)\rightarrow 0
\end{gather}
restricting to $\sA_{E/S}\subset {\sD}iff(E,E)$, and
pushing forward 
by the trace map
$\frac{E\boxtimes E^{\circ}}{E\boxtimes E^{\circ}(-\Delta)}\to \omega_{\Delta}
\cong \omega_{X/S}$, 
yields an exact sequence  
\begin{gather}
0\to \omega_{X/S}\to ^{tr}\sA^{-1}_{E}\xrightarrow{\gamma_E} 
\sA_{E/S}\to 0
\end{gather}
One defines the trace complex $^{tr}\sA^{\bullet}$ 
 by  
$\sA_{E,\pi}$ for $i=0$, 
$^{tr}\sA^{-1}_{E}$ for $i=-1$, 
$\sO_X$ for $i=-2$ and $0$ else, 
with differentials $d^{-1}:=\gamma_E$ and $d^{-2}$  
equal to the relative K\"ahler differential (see \cite{BS}, section 2). 

One has an exact sequence of complexes
\begin{gather}
0 \to \Omega^{\bullet}_{X/S}[2] \to ^{tr}\sA_E^{\bullet} \to 
(T_{X/S} \to T_{\pi})[1] \to 0,
\end{gather}
where $\Omega^{\bullet}_{X/S}$ is the relative de Rham complex of $\pi$.
Taking relative cohomology, one obtains the exact sequence
\begin{gather}
0\to \sO_S\to R^0\pi_*(^{tr}\sA^{\bullet}_{E})\to T_S\to 0. 
\end{gather}
Furthermore $R^0\pi_*(^{tr}{\sA}^{\bullet}_E)$ is
a sheaf of algebras (\cite{BS}, 1.2.3). One denotes by   
$\pi(^{tr}\sA^{\bullet}_E)$ the sheaf on $S$ together with its
algebra structure.

 Finally, let $\sB_i, i=1,2$ be two sheaves of algebras on $S$, with an
exact sequence of sheaves of algebras
\begin{gather}
0 \to \sO_S \to \sB_i \to T_S\to 0.
\end{gather}
One defines $\sB_1 + \sB_2$ by taking the subalgebra of $\sB_1
\oplus \sB_2$, inverse image $\sB_1 \times _{T_S} \sB_2$
of the diagonal embedding $T_S \to
T_S\oplus T_S$, and its push out via the trace map $\sO_S
\oplus \sO_S \to \sO_S$. 

The aim of this note is to prove
\begin{thm} \label{thm} 
Let $D\subset X$ be a divisor, \'etale over
$S$. One has a canonical isomorphism
\begin{gather}
\pi(^{tr}\sA^{\bullet}_E)\cong 
\pi(^{tr}\sA^{\bullet}_{E(-D)})+ 
\sA_{{\rm det}\pi_*(E|_{D})}
\end{gather}
\end{thm}
This is \cite{BS} Theorem 2.3.1, ii). We explain in more details
the proof given there and correct a point in it.
\section{Proof of Theorem  \ref{thm}}

The proof uses the construction of a complex $\sL^{\bullet}$,
together with maps $\sL^{\bullet} \to ^{tr}\sA^{\bullet}_E$ and $
\sL^{\bullet} \to ^{tr}\sA^{\bullet}_{E(-D)} \oplus i_{D*}\sA_{E|_{D}}$
inducing isomorphisms from $R^0\pi_*\sL^{\bullet}$ with the
left and the right hand side of theorem \ref{thm}. We make the
construction of $\sL^{\bullet}$ and the maps explicit, and show
that the induced morphisms are surjective, with the same
(non-vanishing) kernel. 

We first recall the definition 
of the sub-complex
$\sL^{\bullet}\subset ^{tr}\sA^{\bullet}_E$ 
(see \cite{BS}, theorem 2.3.1, ii)):
$\sL^0 \subset ^{tr}\sA^0_E$ consists of the differential 
operators $P$ with $\epsilon (P) \in T_\pi<-D>$, where
$T_\pi<-D>= T_\pi \cap T_X<-D>$ and $T_X<-D>= {\sH}om_{\sO_X}
(\Omega^1_X<D>, \sO_X)$ where 
$\Omega^1_X<D>$ denotes the sheaf of 1-forms with log poles along $D$.   
In particular,
$(d^{-1})^{-1}(\sL^0)\subset ^{tr}\sA_E$ maps to $
i_{D*} \sA_{E|_{D}}$, and $\sL^{-1}\subset (d^{-1})^{-1}(\sL^0)$
is defined as the kernel. Then 
 $\sL^{-2}=\sO_X$. The product struture on
$^{tr}\sA^{\bullet}_E$ is defined in \cite{BS}, 2.1.1.2, and
coincides with the Lie algebra structure on $^{tr}\sA^0_E=
\sA_{E, \pi}$. Since $\sL^{-2}=^{tr}\sA^{-2}_E$,
to see that the product structure stabilizes $\sL^{\bullet}$,
one just has to see that $\sL^0\subset ^{tr}\sA^0_E$ is a  
subalgebra, which is obvious, and that $\sL^0 \times \sL^{-1} \to 
    ^{tr}\sA^{-1}_E$ takes values in $\sL^{-1}$, which is
a consequence of proposition \ref{prop}.

As in section 1, we denote  by $\sA_{E/S}$ the relative Atiyah algebra of
$E$, with symbolic part $T_{X/S}$
and by $\sA_{E,\pi}$ Beilinson's  subalgebra of the global 
Atiyah algebra with symbolic part $T_{\pi}$.  
If $\iota: F\subset E$ is a vector bundle, isomorphic to $E$ away of $D$,
then one has an injection of differential operators
\begin{gather}
{\sD}iff (E, F) \xrightarrow{i} {\sD}iff (E, E)
\end{gather}
induced by $\iota$ on the second argument, and an injection
\begin{gather}
{\sD}iff (E, F) \xrightarrow{j} {\sD}iff (F, F)
\end{gather}
induced by $\iota$ on the first argument.
One has
\begin{defn} \label{EF}
$$\sA_{(E/S, F/S)}:= \sA_{E/S} \cap _i {\sD}iff (E, F) \cong \sA_{F/S} \cap 
_j{\sD}iff (E, F)$$
\end{defn}
Recall $\gamma_{E} : ^{tr}\sA^{-1}_E \to \sA_{E/S}$
denotes the map coming from the filtration by the order of poles of
$\sO_{X\times X}(*\Delta)$ on $^{tr}\sA^{-1}_E$.

One has
\begin{prop} \label{prop}
$$\gamma_E^{-1}( \sA_{E/S, E(-D)/S}) \cong
\gamma_{E(-D)}^{-1}(\sA_{E/S, E(-D)/S}) \cong \sL^{-1}.$$
\end{prop}
\begin{proof}
One considers  
\begin{gather}
\frac{E(-D)\boxtimes E^{\circ}(2\Delta) + E\boxtimes
E^{\circ}}{E\boxtimes E^{\circ}(-\Delta)}\\
= [\frac{E(-D)\boxtimes E^{\circ}(2\Delta)}{E(-D) \boxtimes
E^{\circ}(-\Delta)} \oplus  \frac{E\boxtimes
E^{\circ}}{E\boxtimes E^{\circ}(-\Delta)}]/[\frac{E(-D)\boxtimes
E^{\circ}}{E(-D) \boxtimes
E^{\circ}(-\Delta)}] \notag 
\end{gather}
which, via the natural inclusion to
\begin{gather}
\frac{E\boxtimes E^{\circ}(2\Delta) }{E\boxtimes
E^{\circ}(-\Delta)} 
\end{gather}
is the inverse image  
$\gamma_E^{-1}\Big( {\sD}iff(E, E(-D))\Big)$
(here we abuse of notation, still denoting by $\gamma_E$ the map
coming from the filtration), and via the map coming from the
natural inclusion
\begin{gather}
\frac{E(-D)\boxtimes E^{\circ}(2\Delta)}{E(-D) \boxtimes
E^{\circ}(-\Delta)} \to \frac{E(-D)\boxtimes E^{\circ}(D)
(2\Delta)}{E(-D) \boxtimes
E^{\circ}(D)(-\Delta)}
\end{gather}
and the identification with the first term of the filtration
on $\frac{E(-D)\boxtimes
E^{\circ}(D)(2\Delta)}{E(-D)\boxtimes E^{\circ}(D)(-\Delta)}$
\begin{gather}
\frac{E\boxtimes
E^{\circ}}{E\boxtimes E^{\circ}(-\Delta)} \cong 
\frac{E(-D)\boxtimes
E^{\circ}(D)}{E(-D)\boxtimes E^{\circ}(D)(-\Delta)}
\end{gather}
 is the inverse image
$\gamma_{E(-D)}^{-1}\Big( {\sD}iff(E, E(-D))\Big)$.
\end{proof}
The filtration induced by the order of poles of $\sO_{X\times
X}(*D)$ induces the exact sequences
\begin{gather}
0\to {\sH}om(E, E(-D))\to \sA_{(E/S, E(-D)/S)} \to T_{X/S}(-D)\to
0\\ 
0\to {\sE}nd(E)\to \sA_{E/S} \to T_{X/S}\to 0\\
0\to {\sE}nd(E)\to \sA_{E(-D)/S} \to T_{X/S}\to 0.
\end{gather}
Now, as one has an injection $\sL^{\bullet} \subset
^{tr}\sA^{\bullet}_E$ with cokernel $\sQ$, and again by looking
at the filtration by the order of poles on the sheaf in degree (-1),
one obtains 
\begin{gather}
\sQ\cong {\sE}nd(E)|_{D} [1]
\end{gather}
and 
\begin{thm} \label{thm2.3}
One has an exact sequence
\begin{gather} \notag
0 \to R^0\pi_* ({\sE}nd(E)|_{D}) \to R^0\pi_*(\sL^{\bullet})
\to R^0\pi_*(^{tr}\sA^{\bullet}_E) \to 0.
\end{gather}
\end{thm}
On the other hand, one has an injection $\sL^{\bullet} \subset
^{tr}\sA^{\bullet}_{E(-D)} \oplus {i_D}_*\sA_{E|_D}$ with cokernel $\sP$,
and, as $\sL^{\bullet}$ injects into $^{tr}\sA_{E(-D)}^{\bullet}$, 
one has an exact sequence
\begin{gather}
0\to {i_D}_*\sA_{E|_D}[0] \to \sP \to 
[^{tr}\sA^{\bullet}_{E(-D)}/\sL^{\bullet}]
\to 0,
\end{gather}
where $i_D: D \to X$ is the closed embedding.
We see that the induced filtration on the sheaf in degree (-1)  of
$[^{tr}\sA_{E(-D)}^{\bullet}/\sL^{\bullet}]$ has graded pieces 
$(0, {\sE}nd(E|_D), T_{X/S}|_D)$, whereas the filtration on the
sheaf in degree (0)  has graded pieces
$(0, T_{\pi}/T_{\pi}<-D> \cong T_{X/S}|_D)$. 
This last point comes from the obvious
\begin{lem}
\begin{gather}
\notag \{P\in {\sD}iff(E, E), P(E(-D)) \subset E(-D)\} \cong \notag \\
\notag \{P\in {\sD}iff(E(-mD), E(-mD)), \epsilon(P) \in
{\sE}nd(E)\otimes T<-D>\} \notag 
\end{gather}
for any $m\in \Z$, where $\epsilon$ is the symbol map.
\end{lem}
So 
\begin{lem}
$[^{tr}\sA_{E(-D)}^{\bullet}/\sL^{\bullet}]$ is quasiisomorphic to 
${\sE}nd(E|_D)[1]$.
\end{lem}
The connecting morphism 
$R^{-1}\pi_*[^{tr}\sA_{E(-D)}^{\bullet}/\sL^{\bullet}] \to R^0\pi_*({i_D}_*\sA_{E|_D})[0]$
is just the natural embedding 
$\pi_*({\sE}nd(E|_D)) \to \pi_*({i_D}_*\sA_{E|_D})$ with cokernel
$\pi_* \pi|_{D}^{-1} T_S$. If $D$ is irreducible, 
one has $\pi_* \pi|_{D}^{-1} T_S\cong T_S$, and therefore
\begin{prop}
If $D$ is irreducible, 
one has an exact sequence
\begin{gather} 0 \to R^0\pi_*\sL^{\bullet} \to 
R^0\pi_*[^{tr}\sA^{\bullet}_{E(-D)}]  \oplus
R^0\pi_*[{i_D}_*\sA_{E|_D}]
\to  T_S \to 0\notag
\end{gather}
and the image of $R^0\pi_*\sL^{\bullet}$ is obtained from the
direct sum by taking the pull back under the diagonal embedding
$T_S \to T_S \oplus T_S$.
\end{prop}
On the other hand, still assuming $D$ irreducible, one has the exact sequence
\begin{gather}
0 \to \pi_*{\sE}nd(E|_{D})
\to R^0\pi_*[i_{D*}\sA_{E|_D}] \to T_S\cong \pi_*\pi|_{D}^{-1}T_S \to 0
\end{gather}
and the Atiyah
algebra $\sA_{{\rm det}(\pi_*E|_{D})}$ is the push out of 
$R^0\pi_*[i_{D*}\sA_{E|_D}]$ by the trace map
$\pi_*{\sE}nd(E|_{D}) \to \sO_S$. 

Defining
\begin{gather}
\sK := {\rm Ker}\Big( \sO_S \oplus \pi_*{\sE}nd(E|_{D})
\xrightarrow{{\rm id} \oplus {\rm Tr}} \sO_S\Big) \cong
\pi_*{\sE}nd(E|_D), 
\end{gather}
one thus obtains
\begin{thm}\label{thm2.7}
If $D$ is irreducible, one has an exact sequence
\begin{gather}
0\to \sK \to R^0\pi_*\sL^{\bullet} \to 
\pi_*(^{tr}\sA^{\bullet}_{E(-D)})+
\sA_{{\rm det}\pi_*(E|_D)}\to 0.  \notag
\end{gather}
\end{thm}
It can be easily shown that the embedding  $\pi_*{\sE}nd(E|_D) \subset
R^0\pi_*\sL^{\bullet}$ in theorems \ref{thm2.3} and \ref{thm2.7}
is the same embedding of a subsheaf of ideals.
It finishes the proof of theorem \ref{thm} when 
$D$ is irreducible. In general, since $D$ is \'etale over 
$S$, its irreducible components are disjoint, thus one proves
theorem \ref{thm} by adding one component at a time.

\bibliographystyle{plain}

\begin{thebibliography}{99}
\bibitem{BS} A. Beilinson, V. Schechtman: Determinant Bundles
and Virasoro Algebras, Commun. Math. Phys. {\bf 118} (1988), 651-701.   

\end{thebibliography}
\renewcommand\refname{References}

\end{document}